\documentclass[noae]{amsart}
\usepackage[latin9]{inputenc}
\usepackage{amssymb}
\usepackage{esint}
\begin{document}

\title{On Certain Positive Semidefinite Matrices of Special Functions}

\author{Ruiming Zhang}

\address{College of Science\\
 Northwest A\&F University\\
 Yangling, Shaanxi 712100\\
 P. R. China.}

\email{ruimingzhang@yahoo.com}

\thanks{This work is partially supported by National Natural Science Foundation
of China, grant No. 11371294 and Northwest A\&F University of China.}
\begin{abstract}
Special functions are often defined as a Fourier or Laplace transform
of a positive measure, and the positivity of the measure manifests
as positive definiteness of certain matrices. The purpose of this
expository note is to give a sample of such positive definite matrices
to demonstrate this connection for some well-known special functions
such as Gamma, Beta, hypergeometric, theta, elliptic, zeta and basic
hypergeometric functions. 
\end{abstract}


\maketitle

\section{Introduction}

Recall that for $n\in\mathbb{N}$ and $A=\left(a_{j,k}\right)_{j,k=1}^{n},\ a_{j,k}\in\mathbb{C}$,
$A$ is called positive semidefinite if and only if the quadratic
form ${\displaystyle \sum_{j,k=1}^{n}}a_{j,k}z_{j}\overline{z_{k}}\ge0$
for all $z_{1},\,z_{2},\dots\,z_{n}\in\mathbb{C}$, and it is positive
definite if it is positive semidefinite and ${\displaystyle \sum_{j,k=1}^{n}}a_{j,k}z_{j}\overline{z_{k}}=0$
implies that $z_{1}=\dots=z_{n}=0$. Given two positive semidefinite
$n\times n$ matrices 
\[
A=\left(a_{j,k}\right)_{j,k=1}^{n},\quad B=\left(b_{j,k}\right)_{j,k=1}^{n},\quad a_{j,k},\,b_{j,k}\in\mathbb{C},
\]
it is well-known that the Schur (Hadamard) product $A\circ B=\left(a_{j,k}b_{j,k}\right)_{j,k=1}^{n}$
is also positive semidefinite, and it satisfies \cite{Horn} 
\[
\det\left(A\circ B\right)\ge\det(A)\cdot\det(B).
\]
Since all the minors of a positive semidefinite matrix are nonnegative,
hence a positive semidefinite matrix of special function entries can
yield many inequalities for special functions. 

Given a positive measure $\mu(x)$ on the real line $\mathbb{R}$,
We denote $\mathcal{H}$ as the Hilbert space of $\mu$-square integrable
functions, 
\[
\mathcal{H}=\left\{ f\big|\int_{\mathbb{R}}\left|f(x)\right|^{2}\mu(dx)<\infty\right\} 
\]
endowed with the usual inner product,
\[
<f,g>=\int_{\mathbb{R}}f(x)\overline{g(x)}\mu(dx),\quad f,\,g\in\mathcal{H}.
\]
Let $\left\{ f_{n}(x)\right\} _{n=0}^{N}\subset\mathcal{H}$, where
$N$ may be any nonnegative integer in $\mathbb{N}_{0}$ or equals
to $\infty$, then the finite sections of the Gram matrices \cite{Andrews,Ismail}
\[
G_{n}=\left(<f_{j},\,f_{k}>\right)_{j,k=0}^{n},\quad n=0,\dots,\,N
\]
are positive semidefinite, and they are positive definite if $\left\{ f_{n}(x)\right\} _{n=0}^{N}\subset\mathcal{H}$
are linearly independent.

In this article we shall list some of positive semidefinite matrices
with special function entries. Our method to obtain positive semidefinite
matrices is first to isolate an inner product structure associated
with the special function, then choose a function set to compute the
corresponding Gram matrices, and finally apply Schur product to the
obtained more general positive semidefinite matrices. Even though
the proofs are completely trivial, these positive semidefinite matrices
sometimes may turn out to be very handy. In the following discussion
if any of the formulas below are not specifically referenced, it means
that they can be found in \cite{Erdelyi,Olver}

\section{Main Results}

Recall that the Jacobi $\theta_{3}$-function is defined by
\[
\theta_{3}(z,q)=\sum_{n=-\infty}^{\infty}q^{n^{2}}e^{2\pi inv},\quad z=e^{2\pi iv},\ |q|<1.
\]
For $0<q<1$, define
\[
\mu(x)=\sum_{n=-\infty}^{\infty}q^{n^{2}}\delta\left(x-n\right),
\]
then $\mu(x)$ is a positive measure on $(-\infty,\infty)$. For $n\in\mathbb{N}$,
$c_{1},c_{2},\dots,c_{n}\in\mathbb{C}$ and $v_{j}\in\mathbb{C},\ j=1,\dots,\,n$
we have
\[
\begin{aligned} & \int_{-\infty}^{\infty}\left|\sum_{j=1}^{n}c_{j}e^{2\pi iv_{j}x}\right|^{2}d\mu(x)=\sum_{j,k=0}^{n}c_{j}\overline{c_{k}}\int_{-\infty}^{\infty}e^{2\pi ix(v_{j}-\overline{v_{k}})}d\mu(x)\\
 & =\sum_{j,k=0}^{n}c_{j}\overline{c_{k}}\sum_{\ell=-\infty}^{\infty}q^{\ell^{2}}e^{2\pi i\ell(v_{j}-\overline{v_{k}})}=\sum_{j,k=0}^{n}c_{j}\overline{c_{k}}\theta_{3}\left(e^{2\pi i(v_{j}-\overline{v_{k}})},q\right)\ge0.
\end{aligned}
\]
Thus, for $n\in\mathbb{N}$, $0<q<1$ and $v_{j}\in\mathbb{C},\ j=1,\dots,\,n$,
the matrix
\begin{equation}
\left(\theta_{3}\left(e^{2\pi i(v_{j}-\overline{v_{k}})},\,q\right)\right)_{j,k=1}^{n}\label{eq:m1a}
\end{equation}
is positive semi-definite. For 
\[
v_{j,\ell}\in\mathbb{C},\ 0<q_{\ell}<1,\ 1\le j\le n,\ 1\le\ell\le m,\quad m,\,n\in\mathbb{N},
\]
by taking Schur product of the above matrix we prove that 
\begin{equation}
\left(\prod_{\ell=1}^{m}\theta_{3}\left(e^{2\pi i(v_{j,\ell}-\overline{v_{k,\ell}})},\ q_{\ell}\right)\right)_{j,k=1}^{n}\label{eq:m1b}
\end{equation}
 is also positive semidefinite.

The Jacobi elliptic function $\mbox{dn}(2Kv)$ is defined by
\[
\mbox{dn}(2Kv)=\frac{\pi}{K}\sum_{n\in\mathbb{Z}}\frac{q^{n}}{1+q^{2n}}e^{2n\pi vi},
\]
where 
\[
|q|<1,\ K=\frac{\pi}{2}\theta_{3}^{2}(0,\,q),\ qe^{2\pi\left|\Im(v)\right|}<1.
\]
For $0<q<1$, let
\[
\mu(x)=\frac{\pi}{K}\sum_{n\in\mathbb{Z}}\frac{q^{n}\delta(x-n)}{1+q^{2n}},
\]
is a positive measure on $\mathbb{R}$. For $n\in\mathbb{N}$, $c_{1},\,c_{2},\dots,\,c_{n}\in\mathbb{C}$
and 
\begin{equation}
qe^{4\pi\left|\Im(v_{j})\right|}<1,\ v_{j}\in\mathbb{C},\ j=1,\dots,\,n,\label{eq:15}
\end{equation}
then the quadratic form 
\[
\int_{-\infty}^{\infty}\left|\sum_{j=1}^{n}c_{j}e^{2\pi iv_{j}x}\right|^{2}d\mu(x)=\sum_{j,k=1}^{n}c_{j}\overline{c_{k}}\mbox{dn}\left(2K(v_{j}-\overline{v_{k}})\right)
\]
is nonnegative. Hence, for 
\[
0<q<1,\ K=\frac{\pi}{2}\theta_{3}^{2}(0,\,q),\ qe^{4\pi\left|\Im(v_{j})\right|}<1,
\]
 the matrix 
\begin{equation}
\left(\mbox{dn}\left(2K(v_{j}-\overline{v_{k}})\right)\right)_{j,k=1}^{n}\label{eq:m2a}
\end{equation}
 is positive semidefinite. By taking Schur product we see that 
\begin{equation}
\left(\prod_{\ell=1}^{m}\mbox{dn}\left(2K_{\ell}\left(v_{j,\ell}-\overline{v_{k,\ell}}\right)\right)\right)_{j,k=1}^{n},\quad m,n\in\mathbb{N},\label{eq:m2b}
\end{equation}
is also positive semidefinite where
\[
0<q_{j}<1,\ K_{j}=\frac{\pi}{2}\theta_{3}^{2}(0,\,q_{j}),\ q_{j}e^{4\pi\left|\Im\left(v_{j,\ell}\right)\right|}<1,\ 1\le j\le m.
\]
 The Riemann zeta function is defined as the analytic continuation
of the Dirichlet series \cite{Erdelyi,Olver}
\[
\zeta(s)=\sum_{n=1}^{\infty}\frac{1}{n^{s}},\quad\Re(s)>1.
\]
 For $\Re(s)>0$, it has the following integral representations,

\[
\frac{1}{s-1}-\frac{\zeta(s)}{s}=\int_{1}^{\infty}u^{-s}\frac{\{u\}}{u}dx,\quad\Re(s)>0,
\]
where $0\le\{x\}=x-\left\lfloor x\right\rfloor <1$ is the fractional
part of $x$.

Given $n\in\mathbb{N}$, for $c_{j},\,s_{j}\in\mathbb{C},\ \Re(s_{j})>0,\ j=1,\dots,\,n$,
we have
\[
\int_{1}^{\infty}\left|\sum_{j=1}^{n}\frac{c_{j}}{u^{s_{j}}}\right|^{2}\frac{\{u\}}{u}dx\ge0.
\]
 Then for $n\in\mathbb{N}$, the matrix 
\begin{equation}
\left(\frac{1}{s_{j}+\overline{s_{k}}-1}-\frac{\zeta(s_{j}+\overline{s_{k}})}{s_{j}+\overline{s_{k}}}\right)_{j,k=1}^{n}\label{eq:m3a}
\end{equation}
is positive semidefinite where $\Re(s_{j})>0,\ j=1,\dots,\,n$. 

For $m,\,n\in\mathbb{N}$, by taking Schur product we see the matrix
\begin{equation}
\left(\prod_{\ell=1}^{m}\left\{ \frac{1}{s_{j,\ell}+\overline{s_{k,\ell}}-1}-\frac{\zeta(s_{j,\ell}+\overline{s_{k,\ell}})}{s_{j,\ell}+\overline{s_{k,\ell}}}\right\} \right)_{j,k=1}^{n}\label{eq:m3b}
\end{equation}
 is also positive semidefinite where $\Re(s_{j,\ell})>0,\ j=1,\dots,\,n,\,\ell=1,\dots,\,m$. 

Recall the Euler Gamma function $\Gamma(z)$ is defined as the analytic
continuation of integral,
\[
\Gamma(z)=\int_{0}^{\infty}e^{-x}x^{z-1}dx,\quad\Re(z)>0.
\]
Then for $m,\,n\in\mathbb{N}$, by taking 
\[
f_{k}(x)=x^{z_{k}},\ d\mu(x)=e^{-x}\frac{dx}{x}1_{\{x>0\}},
\]
 then $<f_{j},f_{k}>=\Gamma\left(z_{j}+\overline{z_{k}}\right)$,
hence matrices 
\begin{equation}
\left(\Gamma\left(z_{j}+\overline{z_{k}}\right)\right)_{j,k=1}^{n}\label{eq:m4a}
\end{equation}
and
\begin{equation}
\left(\prod_{\ell=1}^{m}\Gamma\left(z_{j,\ell}+\overline{z_{k,\ell}}\right)\right)_{j,k=1}^{n}\label{eq:m4b}
\end{equation}
are positive semidefinite for $\Re(z_{j}),\ \Re(z_{j,\ell})>0,\ 1\le j\le n,\ 1\le\ell\le m$.

For $\lambda>0$ and $0<\phi<\pi$ we have the following integral
\[
\frac{1}{2\pi}\int_{-\infty}^{\infty}e^{(2\phi-\pi)x}\left|\Gamma(\lambda+ix)\right|^{2}dx=\frac{\Gamma(2\lambda)}{(2\sin\phi)^{2\lambda}},
\]
which is integral of the weight function for the Meixner-Pollaczek
orthogonal polynomials. By taking function sequence $f_{j}(x)=e^{2\phi_{j}x},\ j=1,\dots,\,n$
and $d\mu(x)=e^{-\pi x}\left|\Gamma(\lambda+ix)\right|^{2}dx$ we
see the matrix 
\begin{equation}
\left(\frac{1}{\sin^{\lambda}(\phi_{j}+\phi_{k})}\right)_{j,k=1}^{n}\label{eq:m5a}
\end{equation}
is positive semidefinite where $\lambda>0,\ \pi/2>\phi_{j}>0,\quad j=1,\dots,\,n,\quad n\in\mathbb{N}$.
By taking the Schur product, the matrix
\begin{equation}
\left(\frac{1}{{\displaystyle \prod_{\ell=1}^{m}}\sin^{\lambda_{\ell}}(\phi_{j,\ell}+\phi_{k,\ell})}\right)_{j,k=1}^{n}\label{eq:m5b}
\end{equation}
is also positive semidefinite where 
\[
\frac{\pi}{2}>\phi_{j,\ell}>0,\ \lambda_{j}>0,\quad1\le j\le n,\,1\le\ell\le m,\quad m,\,n\in\mathbb{N}.
\]
Similarly, from the Euler's Beta function
\[
B(p,q)=\int_{0}^{1}x^{p}(1-x)^{q}\frac{dx}{x(1-x)},\quad\Re(p),\,\Re(q)>0,
\]
 we get that for $m,\,n\in\mathbb{N}$ the matrices, by considering
\[
f_{j}(x)=x^{p_{j}}(1-x)^{q_{j}}\ d\mu(x)=\frac{dx}{x(1-x)}1_{\{0<x<1\}},
\]
 we see the matrices
\begin{equation}
\left(B\left(p_{j}+\overline{p_{k}},\ q_{j}+\overline{q_{k}}\right)\right)_{j,k=1}^{n}\label{eq:m6a}
\end{equation}
 and
\begin{equation}
\left(\prod_{\ell=1}^{m}B\left(p_{j,\ell}+\overline{p_{k,\ell}},\ q_{j,\ell}+\overline{q_{k,\ell}}\right)\right)_{j,k=1}^{n},\label{eq:m6b}
\end{equation}
 are positive semidefinite where $\Re(p_{j}),\ \Re(p_{j,\ell}),\ \Re(q_{j}),\ \Re(q_{j,\ell})>0,\ 1\le j\le n,\ 1\le\ell\le m$.

The shifted factorial is defined by \cite{Andrews,Gasper,Ismail}
\[
(a)_{n}=\frac{\Gamma(a+n)}{\Gamma(a)},\ \left(a_{1},\dots,\,a_{m}\right)_{n}=\prod_{k=1}^{m}(a_{k})_{n},
\]
where $a,\,n,\,a_{1},\dots,\ a_{m}\in\mathbb{C}$. Then for $s+1\ge r$
and $a_{1},\dots,\,a_{r},\,b_{1},\dots,\,b_{s}\in\mathbb{C}$, the
hypergeometric series is defined by 
\[
_{r}F_{s}\left(\begin{array}{cc}
\begin{array}{c}
a_{1},\dots,\,a_{r}\\
b_{1},\dots,\,b_{s}
\end{array} & \bigg|z\end{array}\right)=\sum_{n=0}^{\infty}\frac{\left(a_{1},\dots,\,a_{r}\right)_{n}}{\left(1,\,b_{1},\dots,\,b_{s}\right)_{n}}z^{n},
\]
where $z\in\mathbb{C}$ for $s+1>r$ and $|z|<1$ for $s+1=r$. Given
$m,\,n\in\mathbb{N}$ and nonnegative integers $r_{\ell},\:s_{\ell},\ 1\le\ell\le m$,
we assume that $a_{1,\ell},\dots,\,a_{r_{\ell},\ell},\,b_{1,\ell},\dots,\,b_{s_{\ell},\ell},\ 1\le\ell\le m$
are so chosen such that 
\[
\frac{\left(a_{1,\ell},\dots,\,a_{r_{\ell},\ell}\right)_{n}}{\left(b_{1,\ell},\dots,\,b_{s_{\ell},\ell}\right)_{n}}\ge0,\quad1\le\ell\le m,\,n\in\mathbb{N}_{0}.
\]
Then the matrix 
\begin{equation}
\left(\prod_{\ell=1}^{m}{}_{r_{\ell}}F_{s_{\ell}}\left(\begin{array}{cc}
\begin{array}{c}
a_{1,\ell},\dots,\,a_{r_{\ell},\ell}\\
b_{1,\ell},\dots,\,b_{s_{\ell},\ell}
\end{array} & \bigg|z_{j,\ell}\overline{z_{k,\ell}}\end{array}\right)\right)_{j,k=1}^{n},\label{eq:m7}
\end{equation}
is positive semidefinite where for $1\le\ell\le m$ we assume that
$z_{j,\ell}\in\mathbb{C},\ 1\le j\le n$ if $s_{\ell}+1>r_{\ell}$,
and $|z_{j,\ell}|<1,\ 1\le j\le n$ if $s_{\ell}+1=r_{\ell}$. 

From 
\[
\left(1-2^{1-s}\right)\Gamma(s)\zeta(s)=\int_{0}^{\infty}u^{s}\frac{du}{u\left(e^{u}+1\right)},\quad\Re(s)>0,
\]
we see the matrix 
\begin{equation}
\left(\prod_{\ell=1}^{m}\left(1-2^{1-s_{j,\ell}-\overline{s_{k,\ell}}}\right)\Gamma(s_{j,\ell}+\overline{s_{k,\ell}})\zeta(s_{j,\ell}+\overline{s_{k,\ell}})\right)_{j,k=1}^{n}\label{eq:m8}
\end{equation}
is positive semidefinite where $n,\,m\in\mathbb{N}$ and $\Re(s_{j,\ell})>0,\ j=1,\dots n,\,\ell=1,\dots,m,$. 

From 
\[
\left(1-2^{1-s}\right)\Gamma(s+1)\zeta(s)=\int_{0}^{\infty}u^{s}\frac{e^{u}du}{\left(e^{u}+1\right)^{2}},\quad\Re(s)>0
\]
we see that 
\begin{equation}
\left(\prod_{\ell=1}^{m}\left(1-2^{1-s_{j,\ell}-\overline{s_{k,\ell}}}\right)\Gamma(s_{j,\ell}+\overline{s_{k,\ell}}+1)\zeta(s_{j,\ell}+\overline{s_{k,\ell}})\right)_{j,k=1}^{n}\label{eq:m9}
\end{equation}
is positive semidefinite where $\Re(s_{j,\ell})>0,\ j=1,\dots n,\,\ell=1,\dots,m,\ n,m\in\mathbb{N}$. 

Because
\[
\frac{\pi(s)_{p}\zeta(p+s)}{\sin\pi s}=\int_{0}^{\infty}\left((-1)^{p-1}\psi^{(p)}\left(1+x\right)\right)\frac{dx}{x^{s}},\quad\Re(s)\in(0,1),\ p\in\mathbb{N},
\]
 where 
\[
\frac{(-1)^{p-1}\psi^{(p)}\left(1+x\right)}{p!}=\sum_{n=1}^{\infty}\frac{1}{(x+n)^{p+1}}
\]
is positive on $(0,\infty)$. Hence the matrix 
\begin{equation}
\left(\prod_{\ell=1}^{m}\frac{(s_{j,\ell}+\overline{s_{k,\ell}})_{p_{\ell}}\zeta(p_{\ell}+s_{j,\ell}+\overline{s_{k,\ell}})}{\sin\pi(s_{j,\ell}+\overline{s_{k,\ell}})}\right)_{j,k=1}^{n},\label{eq:m10}
\end{equation}
is positive semidefinite where $n,\,m,\,p\in\mathbb{N}$ and $0<\Re(s_{j,\ell})<\frac{1}{2},\ p_{\ell}\in\mathbb{N},\ j=1,\dots,\,n,\,\ell=1,\dots,\,m$.

The Riemann Xi function 
\[
\Xi(z)=-\frac{(1+4z^{2})}{8\pi^{(1+2iz)/4}}\Gamma\left(\frac{1+2iz}{4}\right)\zeta\left(\frac{1+2iz}{2}\right)
\]
is an even entire function of genus $1$. It satisfies 

\[
\Xi(z)=\int_{-\infty}^{\infty}e^{-itz}\phi(t)dt.
\]
 Thus the matrix 
\begin{equation}
\left(\prod_{\ell=1}^{m}\Xi\left(z_{j,\ell}-\overline{z_{k,\ell}}\right)\right)_{j,k=1}^{n}\label{eq:m11}
\end{equation}
 is positive semidefinite where $n,\,m\in\mathbb{N}$ and $z_{j,\ell}\in\mathbb{C},\ j=1,\dots,\,n,\,\ell=1,\dots,\,m$. 

The Hurwitz zeta function $\zeta(s,a)$ is defined as the analytic
continuation of 
\[
\zeta(s,a)=\sum_{n=0}^{\infty}\frac{1}{(n+a)^{s}},\quad\Re(s)>1,\ -a\notin\mathbb{N}_{0.}
\]
For $\Re(s)>0,\ a>0$ and $m,\,n\in\mathbb{N}$, from 

\[
\left(\frac{1}{a^{s}}+\frac{1}{(1+a)^{s}}+\frac{(1+a)^{1-s}}{s-1}-\zeta(s,a)\right)=s\int_{1}^{\infty}\frac{\{x\}dx}{(x+a)^{s+1}},
\]
 we see that
\begin{equation}
\left(\prod_{\ell=1}^{m}\left(\frac{a_{\ell}^{-s_{j,\ell}-\overline{s_{k,\ell}}}+(1+a_{\ell})^{-s_{j,\ell}-\overline{s_{k,\ell}}}-\zeta(s_{j,\ell}+\overline{s_{k,\ell}},a_{\ell})}{s_{j,\ell}+\overline{s_{k,\ell}}}+\frac{(1+a_{\ell})^{1-s_{j,\ell}-\overline{s_{k,\ell}}}(s_{j,\ell}+\overline{s_{k,\ell}})^{-1}}{(s_{j,\ell}+\overline{s_{k,\ell}}-1)}\right)\right)_{j,k=1}^{n},\label{eq:m12}
\end{equation}
is positive semidefinite for $\Re(s_{j,\ell}),\ a_{\ell},\ 1\le j\le n,\,1\le\ell\le m$. 

Since for $\Re(s)>0$ and $a>0$ we have
\[
\frac{\Gamma(s)}{4^{s}}\left(\zeta\left(s,\,\frac{a+1}{4}\right)-\zeta\left(s,\,\frac{a+3}{4}\right)\right)=\int_{0}^{\infty}x^{s}\frac{dx}{2xe^{ax}\cosh x}.
\]
 Then, the matrix
\begin{equation}
\left(\prod_{\ell=1}^{m}\Gamma\left(s_{j,\ell}+\overline{s_{k,\ell}}\right)\left(\zeta\left(s_{j,\ell}+\overline{s_{k,\ell}},\,\frac{a_{\ell}+1}{4}\right)-\zeta\left(s_{j,\ell}+\overline{s_{k,\ell}},\,\frac{a_{\ell}+3}{4}\right)\right)\right)_{j,k=1}^{n}\label{eq:m13}
\end{equation}
is positive semidefinite for $\Re(s_{j,\ell}),\ a_{\ell}>0,\ 1\le j\le n,\,1\le\ell\le m$. 

The Lerch's Transcendent $\Phi(z,s,a)$ is defined as the analytic
continuation of
\[
\Phi(z,s,a)=\sum_{n=0}^{\infty}\frac{z^{n}}{(a+n)^{s}},\quad-z\notin\mathbb{N}_{0},\ |z|<1;\ \Re(s)>1,\ |z|=1.
\]
For $a>0,\ z<1$ and $\Re(s)>0$ it has the following integral representation,
\[
\Gamma(s)\Phi(z,\,s,\,a)=\int_{0}^{\infty}x^{s}\frac{dx}{xe^{ax}\left(1-ze^{-x}\right)}.
\]
Then for $m,\,n\in\mathbb{N}$, the matrix 
\begin{equation}
\left(\prod_{\ell=1}^{m}\Gamma(s_{j,\ell}+\overline{s_{k,\ell}})\Phi(z_{\ell},\,s_{j,\ell}+\overline{s_{k,\ell}},\,a_{\ell})\right)_{j,k=1}^{n}\label{eq:m14}
\end{equation}
is positive semidefinite for $\Re(s_{j,\ell}),\ a_{\ell}>0,\ z_{\ell}<1,\ 1\le j\le n,\,1\le\ell\le m$.

For $q\in\left(0,1\right)$ and $m\in\mathbb{N}$, let \cite{Andrews,Gasper,Ismail}
\[
\left(z;\,q\right)_{\infty}=\prod_{n=0}^{\infty}\left(1-zq^{n}\right),\ (z,\,q)_{n}=\frac{\left(z;\,q\right)_{\infty}}{\left(zq^{n};\,q\right)_{\infty}},\quad z,n\in\mathbb{C},
\]
and 
\[
(z_{1},\,z_{2},\dots,\,z_{m};\,q)_{n}=\prod_{j=1}^{m}(z_{j};\,q)_{n},\ z_{j},\,n\in\mathbb{C}.
\]
For $\alpha_{1},\alpha_{2},\alpha_{3},\alpha_{4}>0$, a weaker form
of the Askey-Wilson beta integral is
\[
\begin{aligned} & \int_{0}^{\pi}\frac{(1-q)^{5}(q;\,q)_{\infty}^{6}\left|(e^{2i\theta};\,q)_{\infty}\right|^{2}d\theta}{\left|(q^{\alpha_{1}}e^{i\theta},\,q^{\alpha_{2}}e^{i\theta},\,q^{\alpha_{3}}e^{i\theta},\,q^{\alpha_{4}}e^{i\theta};\,q)_{\infty}\right|^{2}2\pi}\\
 & =\frac{(1-q)^{2(\alpha_{1}+\alpha_{2}+\alpha_{3}+\alpha_{4})}{\displaystyle \prod_{1\le j<k\le4}}\Gamma_{q}(\alpha_{j}+\alpha_{k})}{\Gamma_{q}(\alpha_{1}+\alpha_{2}+\alpha_{3}+\alpha_{4})}.
\end{aligned}
\]
By taking the function sequence,
\[
u_{k}(\theta)=\frac{1}{\left|(q^{\alpha_{k,1}}e^{i\theta},\,q^{\alpha_{k,2}}e^{i\theta};\,q)_{\infty}\right|^{2}},
\]
and the Schur product, we see the matrix
\begin{equation}
\left(\prod_{\ell=1}^{m}\frac{\Gamma_{q}(\alpha_{j,1,\ell}+\alpha_{k,1,\ell})\Gamma_{q}(\alpha_{j,2,\ell}+\alpha_{k,2,\ell})}{\Gamma_{q}(\alpha_{j,1,\ell}+\alpha_{j,2,\ell}+\alpha_{k,1,\ell}+\alpha_{k,2,\ell})}\right)_{j,k=1}^{n}\label{eq:m15}
\end{equation}
is positive semidefinite for 
\[
m,\,n\in\mathbb{N},\ \alpha_{k,1,\ell},\,\alpha_{k,2,\ell}>0,\quad1\le k\le n,\,1\le\ell\le m.
\]
For $r,\,s\in\mathbb{N}_{0}$, $0<q<1$, and $a_{1},\dots,\,a_{r},\ b_{1},\dots,\,b_{s}\in\mathbb{C}$,
let
\[
\begin{aligned} & _{r}A_{s}^{(\alpha)}(a_{1},\dots,\,a_{r};\,b_{1},\dots,\ b_{s};\,q;\,z)={}_{r}\phi_{s}\left(\begin{array}{cc}
\begin{array}{c}
a_{1},\dots,\,a_{r}\\
b_{1},\dots,\,b_{s}
\end{array} & \bigg|q,\,z\end{array}\right)\\
 & =\sum_{n=0}^{\infty}\frac{(a_{1},\dots,\,a_{r};\,q)_{n}}{\left(b_{1},\dots,b_{s};\,q\right)_{n}}q^{\alpha n^{2}}z^{n}.
\end{aligned}
\]
Then it is clear that for $0<q<1$ and $s+1\ge r$ we have 
\[
\begin{aligned} & _{r}A_{s}^{\left((s+1-r)/2\right)}\left(a_{1},\dots a_{r};\ q,\,b_{1},\dots,\ b_{s}\ ;\,q;\,(-1/\sqrt{q})^{s+1-r}z\right)\\
 & ={}_{r}\phi_{s}\left(\begin{array}{cc}
\begin{array}{c}
a_{1},\dots,\,a_{r}\\
b_{1},\dots,\,b_{s}
\end{array} & \bigg|q,\,z\end{array}\right),
\end{aligned}
\]
where the basic hypergeometric series $_{r}\phi_{s}$ is defined by
\cite{Andrews,Gasper,Ismail}
\[
\begin{aligned}_{r}\phi_{s}\left(\begin{array}{cc}
\begin{array}{c}
a_{1},\dots,\,a_{r}\\
b_{1},\dots,\,b_{s}
\end{array} & \bigg|q,\,z\end{array}\right) & =\sum_{n=0}^{\infty}\frac{(a_{1},\dots,\,a_{r};\,q)_{n}}{(q,b_{1},\dots,\,b_{s};\,q)_{n}}\left(-q^{(n-1)/2}\right)^{n(s+1-r)}z^{n}.\end{aligned}
\]
For $n,\,m\in\mathbb{N}$ and $0<q_{\ell}<1,\ 1\le\ell\le m$, we
assume that $z_{j,\ell}\in\mathbb{C},\ 1\le j\le n$ when $\alpha_{\ell}>0$,
and when $\alpha_{\ell}=0$, we restrict $z_{j,\ell},\ 1\le j\le n$
inside an open disk with certain radius less than $1$ to ensure the
associated series converges. Furthermore, we also assume that 
\[
\frac{\left(a_{1,\ell},\dots,\,a_{r_{\ell},\ell};\,q_{\ell}\right)_{n}}{\left(b_{1,\ell},\dots,\,b_{s_{\ell},\ell};\,q_{\ell}\right)_{n}}\ge0,\quad n\in\mathbb{N}_{0}.
\]
 then the matrix 
\begin{equation}
\left(\prod_{\ell=1}^{m}{}_{r_{\ell}}A_{s_{\ell}}^{(\alpha_{\ell})}\left(\begin{array}{cc}
\begin{array}{c}
a_{1,\ell},\dots,\,a_{r_{\ell},\ell}\\
b_{1,\ell},\dots,\,b_{s_{\ell},\ell}
\end{array} & \bigg|q_{\ell},\,z_{j,\ell}\overline{z_{k,\ell}}\end{array}\right)\right)_{j,k=1}^{n}\label{eq:m16}
\end{equation}
 are positive semidefinite. 

For \cite{Gasper} 
\[
q=e^{2\pi i\sigma},\ p=e^{2\pi i\tau},\quad\Im(\sigma),\,\Im(\tau)>0,
\]
let 
\[
\theta(x;\,p)=(x,\,p/x;\,p)_{\infty},\ \theta(x_{1},\dots,\,x_{m};\,p)=\prod_{k=1}^{m}\theta(x_{k};\,p),\quad m\in\mathbb{N},
\]
\[
\left(a;\,q,\,p\right)_{n}=\begin{cases}
\prod_{k=0}^{n-1}\theta(aq^{k};\,p), & n\in\mathbb{N},\\
1 & n=0\\
1/\prod_{k=0}^{-n-1}\theta(aq^{n+k};\,p),\quad & -n\in\mathbb{N}
\end{cases}
\]
 and
\[
\left(a_{1},\,a_{2},\dots,\,a_{m};\,q,\,p\right)_{n}=\prod_{k=1}^{m}(a_{k};\,q,\,p)_{n},\quad m\in\mathbb{N}.
\]
 For $r,\,s\in\mathbb{N}_{0}$ , the modular series $_{r}E_{s}$ and
$_{r}G_{s}$ are defined by \cite{Gasper}
\[
\begin{aligned} & _{r}E_{s}\left(a_{1},\dots,\,a_{r};\,b_{1},\dots,\,b_{s};\,q,\,p;\,A,\ z\right)\\
 & ={}_{r}E_{s}\left(\begin{array}{cc}
\begin{array}{c}
a_{1},\dots,\,a_{r}\\
b_{1},\dots,\,b_{s}
\end{array} & \bigg|q,\ p;\,A,\,z\end{array}\right)=\sum_{n=0}^{\infty}\frac{\left(a_{1},\,a_{2},\dots,\,a_{r};\,q,\,p\right)_{n}}{\left(q,\,b_{1},\dots,\,a_{s};\,q,\,p\right)_{n}}A_{n}z^{n}
\end{aligned}
\]
 and 
\[
\begin{aligned} & _{r}G_{s}\left(c_{1},\dots,\,c_{r};\,d_{1},\dots,\,d_{s};\,q,\,p;\,B,\ z\right)\\
 & ={}_{r}G_{s}\left(\begin{array}{cc}
\begin{array}{c}
c_{1},\dots,\,c_{r}\\
d_{1},\dots,\,d_{s}
\end{array} & \bigg|q,\ p;\,B,\,z\end{array}\right)=\sum_{n=-\infty}^{\infty}\frac{\left(c_{1},\dots,\,c_{r};\,q,\,p\right)_{n}}{\left(d_{1},\dots,\,d_{s};\,q,\,p\right)_{n}}B_{n}z^{n},
\end{aligned}
\]
where the sequences $A_{n}$ and $B_{n}$ are so-chosen to guarantee
the above series converge in their appropriate domains. 

Given $m\in\mathbb{N},\ 1\le\ell\le m$ we let 
\[
r_{\ell},\,s_{\ell}\in\mathbb{N}_{0},\,q_{\ell}=e^{-2\pi\sigma_{\ell}},\ p_{\ell}=e^{-2\pi\tau_{\ell}},\ \sigma_{\ell},\,\tau_{\ell}>0,
\]
and $a_{1,\ell},\dots,\,a_{r_{\ell},\ell};\,b_{1,\ell},\dots,\,b_{s_{\ell},\ell};\,c_{1,\ell},\dots,\,c_{r_{\ell},\ell};\,d_{1,\ell},\dots,\,d_{s_{\ell},\ell}$
and sequences $\left\{ A_{n,\ell}\right\} _{n=0}^{\infty}$ and $\left\{ B_{n,\ell}\right\} _{n=0}^{\infty}$
to ensure that
\[
\frac{\left(a_{1,\ell},\dots,\,a_{r_{\ell},\ell};\,q_{\ell},\,p_{\ell}\right)_{n}}{\left(q_{\ell},\,b_{1,\ell},\dots,\,a_{s_{\ell},\ell};\,q_{\ell},\,p_{\ell}\right)_{n}}A_{n,\ell}\ge0,\quad n\in\mathbb{N}_{0},
\]
\[
\frac{\left(c_{1,\ell},\dots,\,c_{r_{\ell},\ell};\,q_{\ell},\,p_{\ell}\right)_{n}}{\left(d_{1,\ell},\dots,\,d_{s_{\ell},\ell};\,q_{\ell},\,p_{\ell}\right)_{n}}B_{n,\ell}\ge0,\quad n\in\mathbb{Z}
\]
and the series $_{r_{\ell}}E_{s_{\ell}}\,{}_{r_{\ell}}G_{s_{\ell}}$
are convergent on some symmetric open subsets $S_{\ell}\subset\mathbb{C}$
and $T_{\ell}\subset\mathbb{C}$ with respect to the complex conjugation,
then the matrices
\begin{equation}
\left(\prod_{\ell=1}^{m}{}_{r_{\ell}}E_{s_{\ell}}\left(\begin{array}{cc}
\begin{array}{c}
a_{1,\ell},\dots,\,a_{r_{\ell},\ell}\\
b_{1,\ell},\dots,\,b_{s_{\ell},\ell}
\end{array} & \bigg|q_{\ell},\ p_{\ell};\,A_{\ell},\,z_{j,\ell}\overline{z_{k,\ell}}\end{array}\right)\right)_{j,k=1}^{n}\label{eq:m17}
\end{equation}
 and 
\begin{equation}
\left(\prod_{\ell=1}^{m}{}_{r_{\ell}}G_{s_{\ell}}\left(\begin{array}{cc}
\begin{array}{c}
c_{1,\ell},\dots,\,c_{r_{\ell},\ell}\\
d_{1,\ell},\dots,\,d_{s_{\ell},\ell}
\end{array} & \bigg|q_{\ell},\ p_{\ell};\,B_{\ell},\,w_{j,\ell}\overline{w_{k,\ell}}\end{array}\right)\right)_{j,k=1}^{n}\label{eq:m18}
\end{equation}
are positive semidefinite where $z_{j,\ell}\in S_{\ell},\ 1\le j\le n,\,1\le\ell\le m$
and $w_{j,\ell}\in T_{\ell},\ 1\le j\le n,\,1\le\ell\le m$.

\end{document}